\documentclass[10pt]{amsart}
\usepackage{amssymb}

\newtheorem{theorem}{Theorem}

\newtheorem{corollary}[theorem]{Corollary}

\newtheorem{definition}[theorem]{Definition}
\newtheorem{example}[theorem]{Example}

\newtheorem{lemma}[theorem]{Lemma}

\newtheorem{proposition}[theorem]{Proposition}
\newtheorem{remark}[theorem]{Remark}

\numberwithin{equation}{section}
\newcommand{\removeblock}[1]{\relax}

\begin{document}

\title[SDEs driven by fBm]{Operators associated with a stochastic differential equation driven by fractional Brownian motions}
\author{Fabrice Baudoin, Laure Coutin}
\dedicatory{Laboratoire de Probabilit\'es et
Statistiques\\
Universit\'e Paul Sabatier\\
31062 TOULOUSE Cedex 9 France\\
fbaudoin@cict.fr, coutin@cict.fr}

\begin{abstract}
In this paper, by using a Taylor development type formula, we show
how it is possible to associate differential operators with
stochastic differential equations driven by a fractional Brownian
motion. As an application, we deduce that invariant measures for
such SDEs must satisfy an infinite dimensional system of partial
differential equations.
\end{abstract}
\maketitle

\baselineskip 0.30in

\tableofcontents

\section{Introduction and main result}

A $d$-dimensional fractional Brownian motion with Hurst parameter
$H \in (0,1)$ is a Gaussian process
\[
B_t = (B_t^1,...,B_t^d), \text{ } t \geq 0,
\]
where $B^1,...,B^d$ are $d$ independent centered Gaussian
processes with covariance function
\[
R\left( t,s\right) =\frac{1}{2}\left(
s^{2H}+t^{2H}-|t-s|^{2H}\right).
\]
It can be shown that such a process admits a continuous version
whose paths have $p$ finite variation for $1/p<H$. Let us observe
that for $H= \frac{1}{2}$, $B$ is a standard Brownian motion.

In this paper, we are interested in the study in small times of
stochastic differential equations on $\mathbb{R}^n$
\begin{equation}
\label{SDEyoung} X^{x_0}_t =x_0 + \sum_{i=1}^d \int_0^t V_i
(X^{x_0}_s) dB^i_s
\end{equation}
where the $V_i$'s are $C^{\infty}$-bounded vector fields on
$\mathbb{R}^n$ and $B$ is a $d$ dimensional fractional Brownian
motion with Hurst parameter $H > \frac{1}{3}$.

Let us recall that a smooth vector field $V$ on $\mathbb{R}^{n}$
is simply a smooth map
\[
\begin{array}{llll}
V: & \mathbb{R}^{n} & \rightarrow  & \mathbb{R}^{n} \\
& x & \mapsto  & (v_{1}(x),...,v_{n}(x)).
\end{array}
\]
It defines a differential operator acting on the smooth functions
$f: \mathbb{R}^{n} \rightarrow \mathbb{R}$ as follows:
\[
(Vf) (x)=\sum_{i=1}^n v_i (x) \frac{\partial f}{\partial x_i}.
\]
With this notation, we observe that $V$ is a derivation, that is a
map on $\mathcal{C}^{\infty} (\mathbb{R}^{n} , \mathbb{R} )$,
linear over $\mathbb{R}$, satisfying for $f,g \in
\mathcal{C}^{\infty} (\mathbb{R}^{n} , \mathbb{R} )$,
\[
V(fg)=(Vf)g +f (Vg).
\]

If $H>\frac{1}{2}$, the integrals
\[
\int_0^t V_i (X^{x_0}_s) dB^i_s
\]
are understood in the sense of Young's integration; see \cite{rascanu}, \cite{Yo}
and \cite{Za} . But if $H > \frac{1}{3}$ the integrals that appear
in ($\ref{SDEyoung}$) are understood in the rough paths sense of
Lyons (see \cite{CouQia}). For the convenience of the reader, we
included in an appendix at the end of this paper some results of
rough paths theory that are used in our proofs.

By using \cite{CouQia} and Theorem 6.3.1. pp. 179 of \cite{LyQi},
it is possible to show the existence and the uniqueness of a
process $(X_t^{x_0})_{t \geq 0}$ solving ($\ref{SDEyoung}$).
Observe that from the change of variable formula the process
$(X_t^{x_0})_{t \geq 0}$ is such that for every smooth function $f
: \mathbb{R}^n \rightarrow \mathbb{R}$,
\begin{equation*}
f(X^{x_0}_t) =f(x_0) + \sum_{i=1}^d \int_0^t (V_i f) (X^{x_0}_s)
dB^i_s.
\end{equation*}

We denote by $\mathcal{C}^{\infty}_b( {\mathbb R}^n, {\mathbb R})$
the set of compactly supported smooth functions ${\mathbb R}^n
\rightarrow {\mathbb R}$. If $f \in \mathcal{C}^{\infty}_b(
{\mathbb R}^n, {\mathbb R})$, let us denote
\[
\mathbf{P}_t f (x_0) = \mathbb{E} \left( f (X_t^{x_0} ) \right),
\text{ } t \geq 0,
\]
where $X^{x_0}_t$ is the solution of (\ref{SDEyoung}) at time $t$.

Our main result is the following:
\begin{theorem}
\label{main:theorem} Assume $H > \frac{1}{3}$. There exists a
family $\left( \Gamma^H_k \right)_{k \geq 0}$ of differential
operators such that:
\begin{enumerate}
\item If  $f \in \mathcal{C}^{\infty}_b( {\mathbb R}^n, {\mathbb R})$ and $x \in
{\mathbb R}^n$, then for every $N \geq 0$, when $t \rightarrow 0$
\[
\mathbf{P}_t f(x) =\sum_{k=0}^N t^{2kH} (\Gamma^H_k f)(x)
+o(t^{(2N+1)H});
\]
\item
\[
\Gamma_1^H=\frac{1}{2} \sum_{i=1}^d V_i^2;
\]
\item
\[
\Gamma_2^H=\frac{H}{4} \beta (2H,2H) \sum_{i,j=1}^d
V_i^2V_j^2+\frac{2H-1}{8(4H-1)}\sum_{i,j=1}^d V_i V_j^2 V_i
\]
\[
+\left(\frac{H}{4(4H-1)} -\frac{H}{4} \beta (2H,2H) \right)
\sum_{i,j=1}^d(V_i V_j)^2,
\]
where $\beta(a,b)= \int_0^1 x^{a-1} (1-x)^{b-1} dx$;
\item More generally, $\Gamma_k^H$ is a homogeneous polynomial in the $V_i's$ of
degree $2k$;
\item If $H=\frac{1}{2}$ or $V_i V_j =V_j V_i$ for every $1 \leq i,j \leq d$, then
\[
\Gamma_k^H=\frac{1}{k ! 2^k} \left( \sum_{i=1}^d V_i^2 \right)^k.
\]
\end{enumerate}
\end{theorem}

\begin{remark}
The proof of this Theorem relies on the explicit bound of modulus of continuity
 of the It\^o map (see the appendix). In the book of \cite{LyQi}, this is done only 
for sample paths with $p$ finite variation, with $2 \leq p <3,$. This is the
 only reason why this Theorem is stated for $H > \frac{1}{3}$, but it certainly also holds true
for $H \in \left( \frac{1}{4} , \frac{1}{3} \right]$.
\end{remark}
\begin{remark}
In the case of Brownian motion, some more precise results are available in \cite{Bau}, \cite{Ben2} and \cite{Cast}.
\end{remark}
\section{Commutative case}
In this section, we investigate the simplest case which is the
commutative case. More precisely, we assume throughout the section
that the Lie brackets  $[ V_i ,V_j ]=V_i V_j -V_j V_i=0$, $1 \leq
i , j \leq d$.

For $i=1,...,d$, let us denote by $(e^{tV_i})_{t \in \mathbb{R}}$
the (deterministic) flow associated with the ordinary differential
equation
\[
\frac{dx}{dt}=V_i (x_t).
\]
\begin{proposition}\label{flot}
The flow $\Phi_t$  associated with equation ($\ref{SDEyoung}$) is
given by the formula
\[
\Phi_t = e^{V_1 B^1_t} \circ \cdots \circ e^{V_d B^d_t}.
\]
\end{proposition}
\begin{proof}
Observe first that since the vector fields $V_i$'s are commuting,
the flows $(e^{tV_i})_{t \in \mathbb{R}}$ are also commuting. We
set now for $(x,y) \in \mathbb{R}^n \times \mathbb{R}^d$,
\[
F(x,y)=\left( e^{y_1 V_1} \circ ... \circ e^{y_d V_d} \right) (x).
\]
By applying the change of variable formula, we easily see that the
process $\left( e^{B_t^d V_d} x_0 \right)_{t \geq 0}$ is solution
of the equation
\[
d \left( e^{B_t^d V_d} (x_0) \right)=V_d \left( e^{B_t^d V_d}
(x_0) \right) dB_t^d.
\]
A new application of It\^o's formula shows now that, since $V_d$
and $V_{d-1}$ are commuting,
\[
d \left( e^{B_t^{d-1} V_{d-1} } (e^{B_t^d V_d} x_0) \right)=
\]
\[
V_{d-1} \left( e^{B_t^{d-1} V_{d-1} } (e^{B_t^d V_d} x_0) \right)
 dB_t^{d-1} + V_{d} \left( e^{B_t^{d-1} V_{d-1} } (e^{B_t^d V_d}
x_0) \right) dB_t^{d}.
\]
We deduce hence, by an iterative application of the change of
variable formula that the process $(F(x_0,B_t))_{t \geq 0}$
satisfies
\[
d F(x_0,B_t) = \sum_{i=1}^d V_i ( F(x_0,B_t) ) dB^i_t.
\]
Thus, by pathwise uniqueness for the equation ($\ref{SDEyoung}$),
we conclude that
\[
X_t^{x_0}= F(x_0,B_t), \text{ } t \geq 0.
\]
\end{proof}
\begin{remark}
Observe that the expression
\[
e^{V_1 B^1_t} \circ \cdots \circ e^{V_d B^d_t}
\]
is actually defined for every $H \in (0,1)$. Therefore, in the
commutative case, it makes sense to define solutions of stochastic
differential equations driven by fractional Brownian motions
without restriction on the values of the Hurst parameter $H$, and
without using rough paths theory . For instance, solutions of
one-dimensional equations like
\[
dX_t = \sigma (X_t) dB_t
\]
are well defined for any value of the Hurst parameter, see for instance \cite{nourdin}.

\end{remark}
\begin{corollary}
For any smooth $f : \mathbb{R}^n \rightarrow \mathbb{R}$,
\[
\mathbb{E} \left( f (X^{x_0}_t) \right)=\left( \exp \left(
\frac{1}{2} t^{2H} \sum_{i=1}^d V_i^2 \right) f \right) (x_0).
\]
That is, the function
\[
\varphi (t,x)=\mathbb{E} \left( f (X^{x}_t) \right),
\]
satisfies the partial differential equation
\[
\frac{ \partial \varphi}{\partial t}=H t^{2H-1} \sum_{i=1}^d
(V_i^2 f),
\]
associated with the initial condition
\[
\varphi (0,x) = f(x).
\]
\end{corollary}

\begin{proof}
Observe first that from It\^o's formula for the fractional
Brownian motion, see \cite{nualart-cheridito}
\[
\mathbb{E} \left( e^{V_i B^i_t}f (x_0) \right) = f(x_0) + H
\int_0^t s^{2H-1} \mathbb{E} \left(e^{V_i B^i_s}V_i^2 f
(x_0)\right) ds.
\]
Therefore,
\[
\mathbb{E} \left( e^{V_i B^i_t}f (x_0) \right) = \left( \exp
\left( \frac{1}{2} t^{2H}  V_i^2 \right) f \right) (x_0).
\]
It has been seen that (Proposition \ref{flot})
\[
\Phi_t = e^{V_1 B^1_t} \circ \cdots \circ e^{V_d B^d_t}.
\]
Thus
\[
\mathbb{E} \left( f (X^{x_0}_t) \right)=\left( \exp \left(
\frac{1}{2} t^{2H} \sum_{i=1}^d V_i^2 \right) f \right) (x_0).
\]
\end{proof}
\begin{remark}
So, in the commutative case, there is a Feynman-Kac type formula
for solutions of equations driven by fractional Brownian motions.
It shall be shown later that this type of formula \textbf{only
holds} in the commutative case.
\end{remark}
\begin{example}
Let us consider a one-dimensional stochastic differential equation
of the type
\begin{equation}
\label{oneSDE}
 X^{x_0}_t =x_0+\int_0^t \sigma (X^{x_0}_s) dB_s
\end{equation}
where $\sigma : \mathbb{R} \rightarrow \mathbb{R}$ is a
$C^{\infty}$ bounded function and $B$ is a fractional Brownian
motion with Hurst parameter $H \in (0,1)$. Then, the function
\[
\varphi (t,x)=\mathbb{E} \left( f (X^{x}_t) \right),
\]
satisfies the partial differential equation
\[
\frac{ \partial \varphi}{\partial t}=H t^{2H-1} \sigma^2 (x)
\frac{\partial^2 \varphi}{\partial x^2},
\]
associated with the initial condition
\[
\varphi (0,x) = f(x).
\]
\end{example}

\section{Asymptotic development in small times of $\mathbf{P}_t$}
We now study the generic case of non-commuting vector fields.
Throughout this section we assume $H
>\frac{1}{3}$ and introduce the following notations:
\begin{enumerate}
\item
\[
\Delta^k [0,t]=\{ (t_1,...,t_k) \in [0,t]^k, t_1 < ... < t_k
\};
\]
\item If $I=(i_1,...i_k) \in \{1,...,d\}^k$ is a word with length
$k$,
\begin{equation*}
 \int_{\Delta^k [0,t]}  dB^I= \int_{0 < t_1 < ... < t_k \leq t}
 dB^{i_1}_{t_1}  ...  dB^{i_k}_{t_k}.
\end{equation*}
\end{enumerate}

\begin{theorem}
\label{thme-dl-pt} For  $f \in C^{\infty}_b( {\mathbb R}^n,
{\mathbb R})$ , $x \in {\mathbb R}^n$, and $N \geq 0$, when $t
\rightarrow 0$,
\[
 f(X_t) =f(x) + \sum_{k=1}^N t^{2kH}
\sum_{I=(i_1,...i_{2k} ) } (V_{i_1} ... V_{i_{2k}} f) (x)
\int_{\Delta^{2k} [0,1]}  dB^I 
+o(t^{(2N+1)H}),
\]
and
\[
\mathbf{P}_t f(x) =f(x) + \sum_{k=1}^N t^{2kH}
\sum_{I=(i_1,...i_{2k} ) } (V_{i_1} ... V_{i_{2k}} f) (x)
\mathbb{E}\left(\int_{\Delta^{2k} [0,1]}  dB^I \right)
+o(t^{(2N+1)H}).
\]
\end{theorem}

\begin{proof}
Let us denote by $B^m$ the sequel of linear interpolations of $B$
along the dyadic subdivision of mesh $m$, that is if $t_i^m= i
2^{-m}$ for $i=0,..., 2^m$ then for $t \in [t_i^m, t_{i+1}^m],$
\begin{align*}
B^m_t= B^m_{t_i} + 2^m (t-t_i^m) \left(
B_{t_{i+1}}^m -B^m_{t_i} \right).
\end{align*}
Consider now the equation
\begin{equation}
\label{SDEyoungdiscret} X^{m,x}_t =x + \sum_{i=1}^d \int_0^t
V_i (X^{m,x}_s) dB^{i,m}_s.
\end{equation}
The process $X^{m,x},$ defined in (\ref{SDEyoungdiscret}) has Lipschitz continuous sample paths. Let $p > \frac{1}{H}.$
According to Theorem 5 of \cite{CouQia},
$(X^{m,x}_t,~~t \in [0,1])$ converges to $ (X^{x}_t,~~t \in [0,1])$ in
the distance of $p$ variation (see Appendix \ref{summary-rough-path} (\ref{defdist})  for the definition
of this distance).

Let $f$ be in $\mathcal{C}_b^{\infty}( {\mathbb R}^n, {\mathbb R}).$ Using $2N+1$ times the change of variable formula
\[
f(X^{m,x}_t) =f(x) + \sum_{i=1}^d \int_0^t
(V_i f) (X^{m,x}_s) dB^{i,m}_s.
\]
we obtain
\begin{align}\label{young-ito-approx}
f(X^{m,x}_t)&= f(x) + \sum_{k=1}^{2N+1} \sum_{I= (i_1,...,i_k)} (V_{i_1} ... V_{i_k} f)(x) \int_{\Delta^k[0,t]} dB^{I,m}\\
&+ \sum_{I= (i_1,...,i_{2N+2})} \int_{0<u_1< ...<u_{2N+2}<t}(V_{i_{1}} ...
 V_{i_{2N+2}} f)(X^{m,x}_{u_{1}})  dB^{i_{1}}_{u_{1}}...dB^{i_{2N+2}}_{u_{2N+2}} .\nonumber
\end{align}

By taking the expectation we obtain therefore:
\begin{align}\label{young-ito-approx-expect}
{\mathbb E}\left( f(X^{m,x}_t)\right)&= f(x) + \sum_{k=1}^{2N+1} \sum_{I=(i_1,...,i_k)} (V_{i_1} ... V_{i_k} f)(x){\mathbb E} \left(\int_{\Delta^k[0,t]} dB^{I,m}\right)\\
&+ \sum_{I= (i_1,...,i_{2N+2})}{\mathbb E}\left( \int_{0<u_1<...<u_{2N+2}<t}(V_{i_{1}} ... V_{i_{2N+2}} f)(X^{m,x}_{u_1})  dB^{i_{1}}_{u_{1}}
...dB^{i_{2N+2}}_{u_{2N+2}}\right) .\nonumber
\end{align}
Since $f$ is continuous and bounded the left member of (\ref{young-ito-approx-expect}) converges to ${\mathbf P}_tf(x)$ when $m$ goes to infinity~:
\begin{align}\label{term1-expext}
\lim_{m \rightarrow \infty} {\mathbb E}(f(X^{m,x}_t))= {\mathbf P}_tf(x).
\end{align}

Let now ${\mathbb B}^{m}= (1,{\mathbb B}^{m,1},..., {\mathbb B}^{m,2N+1})$ be the smooth
functional over $B^m $ in the sense of Definition 3.1.1 page 30 of \cite{LyQi}
(see also the Example \ref{smooth-rough-path} of the Appendix \ref{summary-rough-path}).
For $k \leq 2N +1,$ $I=(i_1,...,i_k)$ we have therefore
\begin{align*}
\int_{\Delta^k [0,t]}dB^{m,I}={\mathbb B}^{m,I}.
\end{align*}
 According to Theorem 4 of \cite{CouQia},
${\mathbb B}^{m,2}$ converges, in the distance of $p$ variation $p > \frac{1}{H},$
almost surely and in $\mathbf{L}^2$ to the geometric functional denoted by ${\mathbb B}^{2}$.
For $k \leq 2$, we thus have
\begin{align*}
\int_{\Delta^k [0,t]}dB^{I}={\mathbb B}^{2,I}.
\end{align*}
According to Theorem 3.1.3 of \cite{LyQi} or Theorem \ref{extension-thme} of
the Appendix \ref{summary-rough-path},
the geometric functional  ${\mathbb B}^{2}$ has a unique extension in $C_{0,p}(T^{(2N+1)}({\mathbb R}^n))$ denoted by ${\mathbb B}^{2N+1},$ and for $k \leq 2 N+1,$ $I=(i_1,...,i_k)$ almost surely
\begin{align*}
\int_{\Delta^k [0,t]}dB^{I}={\mathbb B}^{2N+1,I}= \lim_{m \rightarrow \infty}\int_{\Delta^k [0,t]}dB^{m,I}.
\end{align*}
Since $\int_{\Delta^k [0,t]}dB^{m,I}$ belongs to the $k$th Wiener chaos of $B,$
the limit also holds in $L^1$ according to \cite{borell}. That is,
for $k \leq 2N+1,$ $I=(i_1,...,i_k)$
\begin{align*}
{\mathbb E} \left( \int_{\Delta^k [0,t]}dB^{I} \right)={\mathbb E}\left({\mathbb B}^{2N+1,I}\right)
= {\mathbb E}\left(\lim_{m \rightarrow \infty}\int_{\Delta^k [0,t]}dB^{m,I}\right).
\end{align*}

Let us now observe that by symmetry for $k$ an odd integer
\begin{align*}
{\mathbb E} \left(\int_{\Delta^{k}[0,t]} dB^{m,I}\right)=0.
\end{align*}
Using scaling property of fractional Brownian motion, we have for $k$ an even integer
\begin{align*}
{\mathbb E} \left( \int_{\Delta^k [0,t]}dB^{I}\right)
=t^{Hk}{\mathbb E} \left( \int_{\Delta^k [0,1]}dB^{I} \right).
\end{align*}

Define for $k=1,...,N$
\begin{align}\label{def-gamma-k}
\Gamma_k^H&=\sum_{I=(i_1,...,i_{k2})} {\mathbb E}\left(\int_{\Delta^{2k}[0,1]} dB^I\right)
V_{i_{1}}....V_{i_{2k}},\\
&=\lim_{m \rightarrow \infty}\sum_{I=(i_1,...,i_{2k})} {\mathbb E}\left(
\int_{\Delta^{2k}[0,1]} dB^{m,I}\right)V_{i_{1}}....V_{i_{2k}}.\nonumber
\end{align}
Then, the first sum in the right member of (\ref{young-ito-approx-expect})
converges to $\sum_{k=1}^N t^{2kH}(\Gamma_k^Hf)(x),$ that is:
\begin{align} \label{lim_expect-fonct}
\lim_{m \rightarrow \infty}\sum_{k=1}^{2N+1} \sum_{I= (i_1,...,i_k)} (V_{i_1} ... V_{i_k} f)(x)
{\mathbb E} \left(\int_{\Delta^k[0,t]} dB^{I,m}\right)=
\sum_{k=1}^N t^{2kH}(\Gamma_k^H f)(x).
\end{align}
According to Theorem 17 of \cite{coutin-friz-victoir} ,
for any $q \geq 1$ and $p >\frac{1}{H}$ there exists a random variable  $C_{p}$
belonging to $L^q$ such that for any $m,$ $k=1,2,$ $I=(i_1,i_k)$ and $(s,t) \in \Delta_{[0,1]}^2,$
\begin{align*}
|{\mathbb B}^{m, I}_{s,t}| \leq C_{p} |t-s|^{k/p}.
\end{align*}
In what follows $\theta$ and $\kappa$ may vary from lines to lines.
According to Theorem \ref{exist-unicite-ed}, there exist $\theta \geq 1,$ $\kappa \geq 1$ ( depending on $x,$ $V_i,~~i=1,...,d$ and $p),$ such that the geometric functional ${\mathbb Z}^{m,2}$ over $(B^m, X^{m,x})$ is controled by $\omega(t,s)= \kappa (1+ C_{p}^{\theta}) |t-s|$ for any $m\in {\mathbb N}.$
Then, according to Theorem \ref{integration} applied to
\begin{align*}
\alpha(b,x)(v,w)= (v,f(x)v)~~\forall (b,x) \in{\mathbb R}^{d+n},~~~\forall (v,w) \in{\mathbb R}^{d+n}
\end{align*}
there exist $\theta \geq 1,$ $\kappa \geq 1$  ( depending on $x,$ $f,$ $V_i,~~i=1,...,d$ and $p),$
such that the geometric functional ${\mathbb G}^{m,2}$ over $(B^m, f(X^{m,x}))$ is controled by $\omega(t,s)= \kappa (1+ C_{p}^{\theta}) |t-s|,$ for any $m\in {\mathbb N}.$
Therefore, according to Theorem \ref{construct-extension}
there exist $\theta \geq 1,$ $\kappa \geq 1$ ( depending on $x,$ $f,$ $V_i,~~i=1,...d,$ $N$ and $p),$ such that for any $I=(i_1,...,i_{2N+2}),$ $m\in {\mathbb N},$ and $t$
\begin{align}\label{maj-reste}
\left| \int_{ 0 < u_1<...<u_{2N+2}<t} f(X^{m,x}_{u_1})dB^{m,i_1}_{u_1}...dB^{m,i_{2N+2}}_{u_{2N+2}}\right|
\leq \kappa (1+ C_{p}^{\theta}) t^{(2N+2)/p}.
\end{align}
By taking the expectation of each members of (\ref{maj-reste}) and using the fact that
 $C_{p}^{\theta}$ belongs to $L^1,$ we deduce that there exists a constant $\kappa$  depending only on  $x,$ $f,$ $V_i,~~i=1,...d,$ $N$ and $p, $
such that for  $I=(i_1,...,i_{2N+2}),$ $m\in {\mathbb N},$ and $t$
\begin{align}\label{maj-reste-expect}
{\mathbb E}(| \int_{ 0 < u_1<...<u_{2N+2}} f(X^{m,x}_{u_1})dB^{m,i_1}_{u_1}...dB^{m,i_{2N+2}}_{u_{2N+2}<t}| ) \leq \kappa  |t|^{(2N+2)/p}.
\end{align}

Finally, by taking the limit of each term of (\ref{young-ito-approx-expect}) and by using
  (\ref{term1-expext}), (\ref{lim_expect-fonct}) and  (\ref{maj-reste-expect})
we get Theorem \ref{thme-dl-pt}.
\end{proof}

\section{Expectation of iterated integrals of the fractional Brownian motion}
In this section we analyse the coefficients appearing in the differential operators
\[
\Gamma_k^H=\sum_{I=(i_1,...i_{2k} ) } \mathbb{E}\left(\int_{\Delta^{2k} [0,1]}  dB^I \right)
 V_{i_1} ... V_{i_{2k}} ,\text{   } H>\frac{1}{4}.
\]
We shall need the following lemma

\begin{lemma} \label{expectationgaussian}
Let $G=(G_1,...,G_{2k})$ be a centered Gaussian vector. We have
\[
\mathbb{E} \left( G_1 ... G_{2k} \right)=\frac{1}{k! 2^k} \sum_{\sigma \in \mathfrak{S}_{2k}}
\prod_{l=1}^k \mathbb{E} \left( G_{\sigma(2l)} G_{\sigma(2l-1)} \right),
\]
where  $\mathfrak{S}_{2k}$ is the group of the permutations
of the set $\{ 1,...,2k \}$.
\end{lemma}

\begin{theorem}
Assume $H>\frac{1}{2}$. Let $I=(i_1,...,i_{2k})$ be a word, then
\[
\mathbb{E}\left(\int_{\Delta^{2k} [0,1]}  dB^I \right)=
\]
\[
\frac{H^k}{k! 2^k} (2H-1)^k \sum_{\sigma \in \mathfrak{S}_{2k}}
\int_{\Delta^{2k} [0,1]} \prod_{l=1}^k \delta_{i_{\sigma(2l-1)},
i_{\sigma(2l)}} \mid t_{\sigma(2l)} - t_{\sigma(2l-1)}\mid^{2H-2}
dt_1...dt_{2k}
\]
where $\delta_{i,j}$ is the Kronecker's symbol.
\end{theorem}
\begin{proof}
Let us first hint how this formula works with a \textit{heuristic} argument. We have
\[
\mathbb{E}\left(\int_{\Delta^{2k} [0,1]}  dB^I \right)
=\int_{\Delta^{2k} [0,1]}\mathbb{E}\left(  dB^I \right).
\]
But, from the  Lemma \ref{expectationgaussian}
\[
\mathbb{E}\left(  dB^I \right)=\frac{1}{k! 2^k} \sum_{\sigma \in
\mathfrak{S}_{2k}} \prod_{l=1}^k \mathbb{E} \left(
dB^{i_{\sigma(2l)}}_{t_{\sigma(2l)}}
dB^{i_{\sigma(2l-1)}}_{t_{\sigma(2l-1)}} \right).
\]
By using the covariance function of the fractional Brownian
motion, we get therefore
\[
\mathbb{E} \left( dB^{i_{\sigma(2l)}}_{t_{\sigma(2l)}}
dB^{i_{\sigma(2l-1)}}_{t_{\sigma(2l-1)}} \right)
=\delta_{i_{\sigma(2l-1)},
i_{\sigma(2l)}}  H(2H-1) \mid t_{\sigma(2l)} -t_{\sigma(2l-1)} \mid^{2H-2} dt_{\sigma(2l)} dt_{\sigma(2l-1)} ,
\]
which leads to the expected result.

We now turn to the rigorous proof. Let us again denote by $B^m$ the sequel of linear interpolation of $B$
along the dyadic subdivision of mesh $m$ and recall that
\begin{align*}
\mathbb{E} \left( \int_{\Delta^{2k} [0,t]}dB^{I}\right)=
\lim_{m \rightarrow \infty} \mathbb{E} \left( \int_{\Delta^{2k} [0,t]}dB^{m,I} \right).
\end{align*}
>From the Lemma \ref{expectationgaussian}, we have
\[
\mathbb{E} \left( \int_{\Delta^{2k} [0,t]}dB^{m,I} \right)=\frac{1}{k! 2^k}
\sum_{\sigma \in \mathfrak{S}_{2k}}
 \int_{\Delta^{2k} [0,t]} \prod_{l=1}^k  \mathbb{E} \left(
\frac{dB^{m,i_{\sigma(2l)}}}{dt_{\sigma(2l)}} \frac{dB^{m,i_{\sigma(2l-1)}}}{dt_{\sigma(2l-1)}}   \right) dt_1 ... dt_{2k}.
\]
If  $t_{\sigma(2l)} \in [t_i^m,t_{i+1}^m[,$ and $t_{\sigma(2l-1)}  \in [t_j^m,t_{j+1}^m[$,
using the expression of $B^m,$ we get
\begin{align*}
\mathbb{E} \left(
\frac{dB^{m,i_{\sigma(2l)}}}{dt_{\sigma(2l)}} \frac{dB^{m,i_{\sigma(2l-1)}}}{dt_{\sigma(2l-1)}}   \right) &=\delta_{i_{\sigma(2l-1)},
i_{\sigma(2l)}}   2^{2m- 2Hm},~~\mbox{for }~ i=j;\\
&= \delta_{i_{\sigma(2l-1)},
i_{\sigma(2l)}}  2^m {\mathbb E}(\Delta^m_iB \Delta^m_{i-1}B),~~\mbox{for }~~j=i-1;\\
&=\delta_{i_{\sigma(2l-1)},
i_{\sigma(2l)}}  H(2H-1)2^{2m}
\int\int_{[t^m_i,t_{i+1}^m] \times  [t^m_j,t_{j+1}^m]}|x-y|^{2H-2} dxdy  ~~\mbox{for }~~j<i-1.
\end{align*}
Here, $\Delta_i^m B$ is the increment of $B$ between $t_i^m$ and $t_{i+1}^m$ that is,
$\Delta_i^m B=B(t_{i+1}^m) -B(t_i^m).$
By using Cauchy-Schwartz inequality, we have
\begin{align*}
|{\mathbb E}( \Delta_i^m B \Delta_{i-1}B^m)| \leq 2^{-2mH}.
\end{align*}
and the result follows from  the Lebesgue convergence dominated
theorem.

\end{proof}

\begin{lemma}
Let $H >\frac{1}{2}$.
\[
\int_{\Delta^{2} [0,1]} (t_2-t_1)^{2H-2} dt_1 dt_2 =\frac{1}{2H(2H-1)}
\]
\[
\int_{\Delta^{4} [0,1]} (t_4 -t_3)^{2H-2} (t_2-t_1)^{2H-2} dt_1 dt_2 dt_3 dt_4=\frac{\beta(2H,2H)}{4H(2H-1)^2}
\]
\[
\int_{\Delta^{4} [0,1]} (t_4 -t_1)^{2H-2} (t_3-t_2)^{2H-2} dt_1
dt_2 dt_3 dt_4=\frac{1}{8H^2 (2H-1)(4H-1)}
\]
\[
\int_{\Delta^{4} [0,1]} (t_4 -t_2)^{2H-2} (t_3-t_1)^{2H-2} dt_1
dt_2 dt_3
dt_4=\frac{1}{4H(4H-1)(2H-1)^2}-\frac{\beta(2H,2H)}{4H(2H-1)^2}
\]
\end{lemma}
\begin{proof}
The result follows by direct but tedious computations.

First of all, we recall that the beta function
\[
\beta (a,b)=\int_0^1 x^{a-1} (1-x)^{b-1}dx
\]
satisfies:
\[
\beta (a+1,b)=\frac{a}{b} \beta (a,b+1)
\]
\[
\beta (a,b)=\beta(a+1,b)+\beta (a ,b+1).
\]
\begin{enumerate}
\item For the first integral:
\begin{align*}
\int_{\Delta^{2} [0,1]} (t_2-t_1)^{2H-2} dt_1 dt_2 =\int_0^1 \int_0^{t_2} (t_2-t_1)^{2H-2}dt_1dt_2 
 =\frac{1}{2H-1} \int_{0}^1 t_2^{2H-1}dt_2=\frac{1}{2H(2H-1)}.
\end{align*}
\item For the second integral:
\begin{align*}
 & \int_{\Delta^{4} [0,1]} (t_4 -t_3)^{2H-2} (t_2-t_1)^{2H-2} dt_1 dt_2 dt_3 dt_4 \\
=& \int_0^1\int_0^{t_4}\int_0^{t_3}\int_0^{t_2} (t_4 -t_3)^{2H-2} (t_2-t_1)^{2H-2} dt_1 dt_2 dt_3 dt_4 \\
=& \frac{1}{2H(2H-1)} \int_0^1\int_0^{t_4} (t_4 -t_3)^{2H-2} t_3^{2H} dt_3 dt_4 \\
=& \frac{1}{8H^2(2H-1)} \int_0^1 (1-s)^{2H-2}s^{2H}ds \quad (t_3 =st_4 ) \\
=& \frac{\beta(2H,2H)}{4H(2H-1)^2}.
\end{align*}
\item For the third integral
\begin{align*}
 & \int_{\Delta^{4} [0,1]} (t_4 -t_1)^{2H-2} (t_3-t_2)^{2H-2} dt_1
dt_2 dt_3 dt_4 \\
=& \frac{1}{2H-1} \int_0^1\int_0^{t_4}\int_0^{t_3} \left( t_4^{2H-1}-(t_4-t_2)^{2H-1} \right)(t_3-t_2)^{2H-2} dt_2dt_3dt_4
\end{align*}
On one hand,
\begin{align*}
 &\frac{1}{2H-1} \int_0^1\int_0^{t_4}\int_0^{t_3}  t_4^{2H-1}(t_3-t_2)^{2H-2}dt_2dt_3dt_4 \\
=&\frac{1}{(2H-1)^2} \int_0^1 \int_0^{t_4}t_4^{2H-1} t_3^{2H-1}dt_3dt_4 \\
=&\frac{1}{8H^2(2H-1)^2}.
\end{align*}
On the other hand,
\begin{align*}
&\frac{1}{2H-1} \int_0^1\int_0^{t_4}\int_0^{t_3}  (t_4-t_2)^{2H-1} (t_3-t_2)^{2H-2} dt_2dt_3dt_4 \\
=&\frac{1}{2H-1} \int_0^1\int_{t_2}^1\int_{t_3}^{1}  (t_4-t_2)^{2H-1} (t_3-t_2)^{2H-2} dt_4dt_3dt_2 \\
=&\frac{1}{2H(2H-1)}\int_0^1 \int_{t_2}^1 (t_3-t_2)^{2H-2} \left( (1-t_2)^{2H}-(t_3-t_2)^{2H} \right)dt_3dt_2 \\
=&\frac{1}{2H(2H-1)^2} \int_0^1 (1-t_2)^{4H-1} dt_2 -\frac{1}{2H(2H-1)(4H-1)}\int_0^1 (1-t_2)^{4H-1} dt_2.
\end{align*}
By putting things together, we obtain the expected result.
\item For the last integral, the computation follows the same lines so that we do not enter into details.
\end{enumerate}
\end{proof}
As an immediate corollary, we deduce
\begin{corollary}\label{gamma2}
Assume $H>\frac{1}{2}$.
\[
\Gamma_1^H=\frac{1}{2} \sum_{i=1}^d V_i^2
\]

\[
\Gamma_2^H=\frac{H}{4} \beta (2H,2H) \sum_{i,j=1}^d
V_i^2V_j^2+\frac{2H-1}{8(4H-1)}\sum_{i,j=1}^d V_i V_j^2 V_i
+\left(\frac{H}{4(4H-1)} -\frac{H}{4} \beta (2H,2H) \right)
\sum_{i,j=1}^d(V_i V_j)^2.
\]
\end{corollary}
It is interesting to observe that the previous corollary makes
sense and actually also holds true for $H>\frac{1}{4}$.
\begin{theorem}
Assume $H>\frac{1}{4}$, then the conclusions of Corollary \ref{gamma2} are still true.
\end{theorem}
We devote now the end of the section to the proof of this theorem.
The trick is to perform transformations on
\[
{\mathbb E} \left(\int_{\Delta^{2k}[0,t]} dB^{m,I}\right)
\]
before passing to the limit.

Showing that for $H>\frac{1}{4}$,
\[
\Gamma_1^H=\frac{1}{2} \sum_{i=1}^d V_i^2
\]
is easy. Indeed, let $(i_1,i_2) \in \{1,...,d\}^2$ such that
$i_1\neq i_2.$ It is easily seen that $B^{m,i_1}$ and $B^{m,i_2}$
are independent, thus
\begin{align*}
{\mathbb E}\left(\int_{\Delta^{2}[0,1]} dB^{m,I}\right)=0.
\end{align*}
Now, if  $i \in \{1,...,d\}$ and $I=(i,i)$ then
\begin{align*}
{\mathbb E}\left(\int_{\Delta^{2}[0,1]}
dB^{m,I}\right)=\frac{1}{2} {\mathbb E}\left(((
B^{m,i}(1))^2\right)= \frac{1}{2}.
\end{align*}
Since $\Gamma^H_1=\sum_{I=(i,j)} \lim_{m \rightarrow \infty}
{\mathbb E}\left(\int_{\Delta^{2}[0,1]} dB^{m,I}\right)V_jV_i,$ we
get the expected result.

Let now $I=(i,j,k,l) \in \{1,...,d\}^4.$ Recall that $B^m $ is
absolutely continuous with respect to the Lebesgue measure, and
the fourth iterated integral is
\begin{align*}
\int_{\Delta^4 [0,1]}dB^{I,m}=\int_{0 <u_1<u_2<u_3 <u_4<1}
\frac{dB^{m,i}}{du}(u_1)\frac{dB^{m,j}}{du}(u_2)\frac{dB^{m,k}}{du}(u_3)\frac{dB^{m,l}}{du}(u_4)du_1
du_2 du_3 du_4.
\end{align*}
Applying Lemma  \ref{expectationgaussian} yields
\begin{align*}
{\mathbb E}\left( \int_{\Delta^4 [0,1]}dB^{I,m}\right)=
\delta_{i,j}\delta_{l,k}A^{m,1}+\delta_{i,k}\delta_{l,j}A^{m,2}+
\delta_{i,l}\delta_{j,k}A^{m,3} ,
\end{align*}
where
\begin{align*}
A^{m,1}= \int_0^1 du_4 \int_0^{u_4} du_3 {\mathbb
E}(\frac{dB^{m}}{du}(u_3)\frac{dB^{m}}{du}(u_4)) \int_0^{u_3}du_2
\int_0^{u_2} du_1 {\mathbb
E}(\frac{dB^{m}}{du}(u_2)\frac{dB^{m}}{du}(u_1));
\end{align*}
\begin{align*}
A^{m,2}= \int_0^1 du_4 \int_0^{u_4} du_2 {\mathbb
E}(\frac{dB^{m}}{du}(u_2)\frac{dB^{m}}{du}(u_4))
\int_{u_2}^{u_4}du_3 \int_0^{u_2} du_1{\mathbb
E}(\frac{dB^{m}}{du}(u_3)\frac{dB^{m}}{du}(u_1));
\end{align*}
\begin{align*}
A^{m,3}= \int_0^1 du_4 \int_0^{u_4} du_1 {\mathbb
E}(\frac{dB^{m}}{du}(u_4)\frac{dB^{m}}{du}(u_1))
\int_{u_1}^{u_4}du_3 \int_{u_1}^{u_3} du_2{\mathbb
E}(\frac{dB^{m}}{du}(u_3)\frac{dB^{m}}{du}(u_2)).
\end{align*}
Therefore, we have to prove that
\begin{align}\label{eq:lim-un}
\lim_{m \rightarrow \infty} A^{m,1}= \frac{H}{4}\beta(2H,2H);
\end{align}
\begin{align}\label{eq:lim-deux}
\lim_{m \rightarrow \infty} A^{m,2}=\frac{H}{4(4H-1)} -\frac{H}{4}
\beta (2H,2H) ;
\end{align}
and
\begin{align}\label{eq:lim-trois}
\lim_{m \rightarrow \infty} A^{m,3}= \frac{2H-1}{8(4H-1)}.
\end{align}

The proof relies on several Lemma.

\begin{lemma}\label{lem-maj-holder}
For any $\alpha <H,$ the random variable
\begin{align*}
\sup_m \|B^m -B\|_{\alpha}= \sup_{m \in {\mathbb N}} \sup_{(s,t)
\in [0,1]^2, s<t} \frac{|[B(t)-B^m(t)]- [B(s)
-B^m(s)]|}{|t-s|^{\alpha}}
\end{align*}
has a finite exponential moment.
\end{lemma}
\begin{proof}
According to Kolmogorov lemma, the random variable $\sup_{(s,t)
\in [0,1]^2, s<t} \frac{|B(t)- B(s)| }{|t-s|^{\alpha}}$ is finite.
Since $B$ is Gaussian process, according to Theorem 1.2.3 of
\cite{fernique} , it has a finite exponential moment. Then the
Lemma is a consequence of results on linear interpolation of
H\"older functions.
\end{proof}
\begin{lemma}\label{cv-un}
For $H > \frac{1}{4},$
\begin{align*}
\lim_{m \rightarrow \infty} A^{m,1}= \frac{H}{4}\beta(2H,2H).
\end{align*}
\end{lemma}
\begin{proof}
Integrating with respect to $u_1$ and $u_2$ in the expression of
$A^{m,1}$ leads to
\begin{align*}
A^{m,1}= \frac{1}{2} \int_0^1 du_4 \int_0^{u_3} du_3
\frac{d^2}{du_3 d u_4} {\mathbb E}(B^{m}(u_3)B^{m}(u_4)) {\mathbb
E}( B^m(u_3)^2).
\end{align*}
Using Fubini's theorem, we integrate first with respect to $u_4$:
\begin{align*}
A^{m,1}&= \frac{1}{2} \int_0^1 du_3  {\mathbb E}([B^m(1) -B^{m}(u_3)]\frac{d}{du }B^{m}(u_3)) {\mathbb E}( B^m(u_3)^2);\\
&= \frac{1}{2}\left[ \int_0^1 du_3 \frac{d}{du_3 } {\mathbb
E}(B^{m}(u_3)B^m(1) ){\mathbb E}( B^m(u_3)^2) -\frac{1}{4}
{\mathbb E}(B^m(1)^2)^2 \right].
\end{align*}
Now, using the expression of $B^m$ we obtain that for $u \in
[t_i^m,t^m_{i+1}[,$
\begin{align*}
 \frac{d}{du } {\mathbb E}(B^{m}(u)B^m(1))&= 2^m{\mathbb E}( \Delta_i^m B B(1))\\
&=\frac{2^m}{2}[|t_{i+1}^m|^{2H} - |1-t_{i+1}^m|^{2H} -(|t_{i}^m|^{2H} - |1-t_{i}^m|^{2H})]\\
&=H2^m \int_{t_i^m}^{t_{i+1}^m}[ r^{2H-1}+ |1- r|^{2H-1}] dr .
\end{align*}
Then, using Fubini's theorem, we get
\begin{align*}
\int_0^1 du {\mathbb E}(B^m(1)\frac{d}{du}B^{m}(u) ) {\mathbb
E}(B^{m}(u)^2)=H \int_0^1 [ r^{2H-1}+ |1- r|^{2H-1}]a ^m(r) dr
\end{align*}
where
\begin{align*}
a^m(r)=  \sum_{i=0}^{2^m-1} 2^m \int_{t_i^m}^{t_{i+1}^m}  {\mathbb
E}(B^{m}(u)^2)du {\mathbf 1}_{[t_i^m, t_{i+1}^m[}(r).
\end{align*}
For all $r\in [0,1],$  $a^m(r)$ converges to ${\mathbb E}(
B(r)^2)=r^{2H}.$  Following Lemma \ref{lem-maj-holder}, it is
bounded uniformly on $m$ and $r.$  Then, using dominated
convergence Lebesgue theorem we obtain
\begin{align*}
\lim_{m \rightarrow +\infty} \int_0^1 du {\mathbb
E}(B^m(1)\frac{d}{du }B^{m}(u) ) {\mathbb E}(B^{m}(u)^2)&=H
\int_0^1 [ r^{2H-1}+ |1- r|^{2H-1}]r^{2H}dr\\
&= \frac{1}{4} + \frac{H}{2} \beta(2H, 2H)
\end{align*}
because $\beta (2H,2H+1)=\frac{1}{2} \beta(2H,2H)$.
We conclude that
\begin{align*}
\lim_{m \rightarrow \infty} A^{m,1} = \frac{H}{4}\beta(2H, 2H).
\end{align*}
\end{proof}

\begin{lemma}
For $H > \frac{1}{4},$
\begin{align*}
\lim_{m \rightarrow \infty} A^{m,3}= \frac{2H-1}{8(4H-1)}.
\end{align*}
\end{lemma}
\begin{proof}
Recall that

\begin{align*}
A^{m,3}= \int_0^1 du_4 \int_0^{u_4} du_1 {\mathbb
E}(\frac{dB^{m}}{du}(u_4)\frac{dB^{m}}{du}(u_1))
\int_{u_1}^{u_4}du_3 \int_{u_1}^{u_3} du_2{\mathbb
E}(\frac{dB^{m}}{du}(u_3)\frac{dB^{m}}{du}(u_2)).
\end{align*}

Integrating with respect to $u_2$ and $u_3$ we obtain
\begin{align*}
A^{m,3}= \frac{1}{2}\int_0^1 du_4 \int_0^{u_4} du_1 {\mathbb
E}(\frac{dB^{m}}{du}(u_4)\frac{dB^{m}}{du}(u_1)) {\mathbb
E}([B^m(u_4)-B^m(u_1)]^2).
\end{align*}
We introduce the indices $i,j$ such that $u_4 \in
[t_i^m,t_{i+1}^m[,$ and $u_1 \in [t_j^m,t_{j+1}^m[.$ Recall that
\begin{align*}
{\mathbb E}(\frac{dB^{m}}{du}(u_1)\frac{dB^{m}}{du}(u_4))&= 2^{2m- 2H},~~\mbox{for }~ i=j;\\
&=2^{2m}{\mathbb E}(\Delta_i^mB \Delta_{i-1}^m B)~~\mbox{ for }~~j=i-1;\\
&=H(2H-1)2^{2m} \int\int_{[t^m_i,t_{i+1}^m] \times
[t^m_j,t_{j+1}^m]}|x-y|^{2H-2} dxdy  ~~\mbox{for }~~j<i-1.
\end{align*}
Then, we split $A^{m,3}$ into three parts
\begin{align*}
A^{m,3}=\frac{1}{2} \left( A^{3,1,m} + A^{3,2,m}+ A^{3,3,m}
\right)
\end{align*}
where
\begin{align*}
A^{3,1,m}=\sum_{i=0}^{2^m-1} 2^{4m- 2Hm} \int_{t_i^m}^{t_{i+1}^m
}du \int_{t_i^m}^u dv (u-v)^2 {\mathbb E}( \Delta_i^m B^2 );
\end{align*}
\begin{align*}
A^{3,2,m}=\sum_{i=0}^{2^m-1} 2^{2m} {\mathbb E}( \Delta_i^mB
\Delta_{i-1}^mB)\int_{t_i^m}^{t_{i+1}^m }du \int_{t_{i-1}^m}^{t_i^m}dv
{\mathbb E}( [B^m(u)-B^m(v)]^2 );
\end{align*}
and
\begin{align*}
A^{3,3,m}=H(2H-1)\sum_{i=2}^{2^m-1} \sum_{j=0}^{i-2} 2^{2m}
\int_{t_i^m}^{t_{i+1}^m }dx  \int_{t_j^m}^{t_{j+1}^m} dy
(x-y)^{2H-2}\int_{t_i^m}^{t_{i+1}^m }  \int_{t_j^m}^{t_{j+1}^m}
dudv
 {\mathbb E}( [ B^m (u) -B^m(v)]^2 ) .
\end{align*}
We have to prove that for $k=1,2$
\begin{align}\label{eq:lim-a3-1}
\lim_{m \rightarrow \infty}A^{3,k,m}=0
\end{align}
and

\begin{align}\label{eq:lim-a3-2}
\lim_{m \rightarrow \infty}A^{3,3,m}=\frac{2H-1}{4(4H-1)}.
\end{align}
First we prove (\ref{eq:lim-a3-1}). Integrating with respect to
$u$ and $v$ yields
\begin{align*}
A^{3,1,m}=\frac{1}{12}\sum_{i=0}^{2^m-1} 2^{- 4Hm } .
\end{align*}

Since $H > \frac{1}{4},$ we deduce that $\lim_{m \rightarrow \infty}A^{3,1,m}=0.$\\
According to Lemma \ref{lem-maj-holder}, we have for $\alpha<H$
\begin{align*}
|A^{3,1,m}|\leq \frac{2}{(2\alpha+1)(2\alpha+2)} 2^{(1 -2H -2
\alpha) m } {\mathbb E}(\sup_m \sup_{(s,t) \in [0,1]^2,
s<t}\frac{|B(s) -B(t)|^2}{|t-s|^{2\alpha}}) .
\end{align*}
Since $H > \frac{1}{4},$ we deduce that $\lim_{m \rightarrow \infty}A^{3,2,m}=0.$\\

Now, we prove (\ref{eq:lim-a3-2}). Indeed using Fubini's theorem,
we have
\begin{align*}
A^{3,3,m}= H(2H-1) \int_0^1 dx \int_0^x (x-y)^{2H-2}a^{3,3,m}(x,y)
dy,
\end{align*}
where for $0 \leq j <i-1 ,$  and $(x,y) \in [t_i^m,t_{i+1}^m[
\times [t_j^m,t_{j+1}^m[$
\begin{align*}
a^{3,3,m}(x,y)= 2^{2m} \int_{t_i^m}^{t_{i+1}^m}
du\int_{t_j^m}^{t_{j+1}^m}dv {\mathbb E}(  (B^m(v) -B^m(u))^2)
\end{align*}
and $a^{3,3,m}(x,y)=0$ elsewhere. For all $x,y$  $a^{2,3,m}(x,y)$
converges almost surely to ${\mathbb E}([B(x)-B(y)]^2).$ Moreover,
since $|t_{i+1}^m -t_j^m|\leq 2|x-y|,$ $a^{3,3,m}$ is bounded by
\begin{align*}
|a^{3,3,m}(x,y)|\leq 2{\mathbb E} (\sup_m \sup_{s<t}\frac{
|B(s) -B(t)|}{|t-s|^{\alpha}})^2)|2(x-y)|^{2\alpha}.
\end{align*}
Using  Lemma \ref{lem-maj-holder} and dominated Lebesgue
convergence theorem, we take the limit when $m$ goes to infinity
for $H > \frac{1}{4}$:
\begin{align*}
\lim_{m \rightarrow \infty} A^{3,3,m}&= H(2H-1) \int_0^1 dx\int_0^x (x-y)^{2H-2}{\mathbb E}([B(x) -B(y)]^2)   dy\\
&=  \frac{2H-1}{4(4H-1)} .
\end{align*}
\end{proof}
Finally,
\begin{lemma}
For $H > \frac{1}{4}$,
\[
\lim_{m \rightarrow \infty} A^{m,2}=\frac{H}{4(4H-1)} -\frac{H}{4}
\beta (2H,2H)
\]
\end{lemma}
\begin{proof}
We know that if the vector fields commute then
\[
\Gamma_2^H=\frac{1}{8} \left( \sum_{i=1}^d V_i^2 \right)^2.
\]
Therefore
\[
\lim_{m \rightarrow + \infty} A^{m,1}
+A^{m,2}+A^{m,3}=\frac{1}{8}.
\]
\end{proof}

\section{Application to the study of invariant measures}

Consider the stochastic differential equation on $\mathbb{R}^n$
\begin{equation}
\label{SDEyoung2} dX_t = V_i (X^{x_0}_t) dB^i_t
\end{equation}
where the $V_i$'s are $C^{\infty}$-bounded vector fields on
$\mathbb{R}^n$ and $B$ is a $d$ dimensional fractional Brownian
motion with Hurst parameter $H > \frac{1}{3}$.

\begin{proposition}
Assume that $\mu$ is a probability measure on $\mathbb{R}^n$ that
is invariant by (\ref{SDEyoung2}). In the sense of distributions,
we have for every $k \geq 1$,
\[
(\Gamma_k^H)^{\ast} \mu =0,
\]
where $(\Gamma_k^H)^{\ast}$ is the formal adjoint of $\Gamma_k^H$.
\end{proposition}

\begin{proof}
Let $ f \in \mathcal{C}_b^{\infty} (\mathbb{R}^n,\mathbb{R})$. We
have for every $N \geq 0$, when $t \rightarrow 0$
\[
\mathbb{E} \left( \int_{\mathbb{R}^n} f(X_t^x) \mu (dx) \right)
=\sum_{k=0}^N t^{2kH} \int_{\mathbb{R}^n} (\Gamma^H_k f)(x) \mu
(dx) +o(t^{(2N+1)H}).
\]
But since $\mu$ is invariant,
\[
\mathbb{E} \left( \int_{\mathbb{R}^n} f(X_t^x) \mu (dx)
\right)=\int_{\mathbb{R}^n} f(x) \mu (dx).
\]
The result follows therefore.
\end{proof}

\section{ Appendix~: On rough path theory}\label{summary-rough-path}

We recall some definitions and  results  on rough path theory, see
\cite{LyQi}. Indeed, we precise how all the constants appearing in
the continuity Theorem of the It\^o map depend on the vectors
fields and the control of the paths.
\subsection{Basic definitions and properties}
We work on $V={\mathbb R}^d$ endowed with the Euclidean norm. The
tensor product is $V^{\otimes k}= V \otimes ...\otimes V$ (of $k$
copies of $V)$ endowed with a norm $|.|_k$ compatible with the
tensor product that is
\begin{align*}
| \xi \otimes \eta|_{k+l} \leq |\xi|_k | \eta|_l,~~~~ \forall \xi \in V^{\otimes k},~~~\forall \eta \in V^{\otimes l}.
\end{align*}
For each $n\in {\mathbb N},$ the truncated tensor algebra $T^{(n)}(V)$ is
\begin{align*}
T^{(n)} (V)= \sum_{k=0}^n V^{\otimes k},~~~V^{\otimes 0}= {\mathbb R}.
\end{align*}
Its multiplication is
\begin{align*}
( \xi \otimes \eta)^k= \sum_{j=0}^k \xi^j \otimes \eta^{k-j},~~~k=0,...n,~~\forall \xi,~ \eta \in T^{(n)}(V).
\end{align*}
The norm $|.|$ on $T^{(n)}(V)$ is defined by
\begin{align*}
|\xi|= \sum_{i=0}^n  | \xi^i|_i,~~~\mbox{ if } \xi=(\xi^0,...,\xi^n).
\end{align*}
The pair $(T^{(n)}(V),~|.|)$ is a tensor algebra with identity
element $(1,0,...,0)$ and for $\xi,\eta \in T^{(n)}(V),~~| \xi
\otimes \eta |\leq |\xi| |\eta|.$

The tensor algebra $T^{(\infty)} (V) $ is
\begin{align*}
T^{( \infty)} (V)= \sum_{k=0}^{\infty} V^{\otimes k},~~~V^{\otimes 0}= {\mathbb R}.
\end{align*}
We use $\Delta^2_{[0,T]} $ to denote the simplex $\{(s,t),~~ 0 \leq s \leq t \leq T\}.$ Recall that a control $\omega$ is a continuous, supper additive function on $\Delta^2_{[0,T]}$ with values in $[0, + \infty[$ such that $\omega(t,t)=0.$
Therefore
\begin{align*}
\omega(s,t) + \omega(t,u) \leq \omega(s,u)~~~\forall (s,t),~~(t,u) \in \Delta^2_{[0,T]}.
\end{align*}
\begin{definition}
A continuous map $X$ from the simplex $\Delta^2_{[0,T]}$ into a truncated tensor algebra $T^{(n)}(V),$ and written as
$X_{s,t}= (X^0_{s,t}, .., X^n_{s,t}) $ with $X^k_{s,t} \in V^{\otimes k}$ for any $( s,t) \in \Delta^2_{[0,T]},~~k=1,...,n$ is called a multiplicative functional of degree $n$ $( n \in {\mathbb N}^*)$ if
\begin{align}\label{chenrule}
&X^0_{s,t}=1 \nonumber \\
&X_{s,t} \otimes X_{t,u}= X_{s,u},~~~\forall (s,t) ,~~(t,u) \in \Delta_{[0,T]}^2,
\end{align}
where the tensor product $\otimes $ is taken in $T^{(n)}(V).$
\end{definition}
Equality (\ref{chenrule}) is called Chen identity, although it appears long before Chen's fundamental works in which a connection is made from iterated path integrals along smooth paths to a class of differential forms on a space of loops on manifold.
\begin{example}\label{smooth-rough-path}
Let $x~: ~ [0,T]\rightarrow V$ be a continuous path. If $x$ is Lipschitz path, then we may build a sequence of iterated integral
$X^k_{s,t}= \int_{s < t_1<...< t_k<t} dx_{t_1} \otimes...\otimes dx_{t_k}.$
In this case identity ( \ref{chenrule}) is equivalent to the additive property of iterated path integrals over differents domains.
\end{example}
\begin{definition}
 Let $p \geq 1$ be a constant. We say that a map $X~:~ [0,T] \rightarrow T^{(n)}(V)$ possesses finite $p$ variation of
\begin{align*}
|X^i_{s,t}| \leq \omega^{i/p}(s,t), ~~\forall i=1,...,n,~~~\forall (s,t) \in \Delta^2_{[0,T]}
\end{align*}
for some control $\omega.$
\end{definition}
\begin{definition}
A multiplicative functional with finite $p$ variation on $T^{([p])}(V)$ is called a rough path (of roughness $p).$ We say that a rough path (of roughness $p)$ $X$ in $T^{([p])}(V)$ is controled by $\omega$ if
\begin{align*}
|X^i_{s,t} | \leq \omega(s,t)^{i/p},~~~ \forall i=1,...,[p],~~~\mbox{ and } \forall ~ (s,t) \in \Delta^2_{[0,T]}.
\end{align*}
\end{definition}

The set of all rough path with roughness $p$  in $T^{([p])}(V)$ will be denoted by $\Omega_p(V).$
\begin{definition}
A smooth rough path $X$ is an element of $\Omega_p(V)$ such that
there exists a   Lipschitz $x~: ~ [0,T]\rightarrow V$  Lipschitz
such that
\begin{align*}
X_{s,t}^k = \int_{s < t_1<...< t_k<t} dx_{t_1} \otimes...\otimes dx_{t_k},~~k=1,...,[p],~~\forall (s,t) \in \Delta^2_{[0,T]}.
\end{align*}
\end{definition}
\subsection{Extension of rough path}
The following theorem shows that the higher (than $[p]$) order
terms $X^k$ $(k > [p])$ are determined uniquely by $X^i~~(i \leq
[p])$ among all possible extensions to a multiplicative functional
wich posses finite $p-$ variation.
\begin{theorem} Let $ p \geq 1,$ and let $X~: ~ \Delta^2_{[0,T]} \rightarrow T^{(n)} (V)$ be a multiplicative functional with finite $p$ variation so that
\begin{align*}
|X_{s,t}^i| \leq \omega(s,t)^{i/p},~~~\forall i=1,..,n ~~\mbox{ and } \forall (s,t) \in \Delta^2_{[0,T]},
\end{align*}
for some control $\omega.$ If $n \geq [p],$ then we may uniquely extend $X$ to be a multiplicative functional in $T^{(\infty)}(V)$ with finite $p$ variation. Moreover, if $\omega$ is a control such that
\begin{align} \label{3.18}
|X^i_{s,t}| \leq  \frac{\omega(s,t)^{i/p}}{ \beta ( \frac{i}{p}) !},~~\forall i=1,...,[p],~~\mbox{ and } \forall (s,t) \in \Delta^2_{[0,T]},
\end{align}
where $\beta$ is a constant such that $$\beta \geq 2 p^2 [ 1+ \sum_{r=3}^{\infty}( \frac{2}{r-2})^{\frac{[p]+1}{p}}],$$
then (\ref{3.18}) remains true for all $i > [p].$
\end{theorem}
The extension of a rough path $X$ to a higher-order multiplicative functional is continuous in $p-$ variation distance.
\begin{theorem}\label{extension-thme} Theorem 3.1.3 p 39 of \cite{LyQi}\\
Let $X$ and $Y$ be two rough path of roughness $p$ and let $\beta$ be a constant such that $\beta \geq 2 p^2 [ 1+ \sum_{r=3}^{\infty}( \frac{2}{r-2})^{\frac{[p]+1}{p}}].$

If $\omega$ is a control such that
\begin{align} \label{3.35}
|X^i_{s,t}|,~~|Y_{s,t}^i | \leq \frac{\omega(s,t)^{i/p}}{ \beta ( \frac{i}{p} )!},~~\forall i=1,...,[p],~~\mbox{ and } \forall (s,t) \in \Delta^2_{[0,T]},
\end{align}
\begin{align} \label{3.36}
|X^i_{s,t}-Y_{s,t}^i | \leq \varepsilon  \frac{\omega(s,t)^{i/p}}{ \beta ( \frac{i}{p}) !},~~\forall i=1,...,[p],~~\mbox{ and } \forall (s,t) \in \Delta^2_{[0,T]},
\end{align}
then (\ref{3.35}) and (\ref{3.36}) hold for  all $i.$
\end{theorem}
\subsubsection{Almost rough path}
In this section we give a method of constructing rough paths.
\begin{definition}
Let $p \geq 1$ be a constant. A function $X~:~\Delta^2\rightarrow T^{([p])}(V)$ is called an almost rough path ( of roughness $p)$ if it is of finite $p-$ variation, $X^0_{s,t}=1$ and for some control $\omega$ and some constant $\theta>1,$
\begin{align*}
|(X_{s,t}\otimes X_{t,u})^i -X_{s,u}| \leq \omega(s,u)^{\theta}
\end{align*}
for all $(s,t) , ~(t,u) \in \Delta^2_{[0,T]}$ and $i=1,..., [p].$
\end{definition}
The following theorem justifies the name of almost rough path.
 \begin{theorem}\label{construct-extension}
If $X~: \Delta^2 \rightarrow T^{([p])}(V)$ is an almost rough path
of roughness $p,$ controlled by $\omega$ and $\theta$ then there
exists an unique rough path $\hat{X}$ (with roughness $p)$ in
$T^{([p])} (V)$ such that
\begin{align*}
| \hat{X^i}_{s,t} - X^i_{s,t}| \leq  K_i\omega(s,t)^{\theta},~~\forall 1 \leq i \leq [p],~~ \forall (s,t) \in \Delta^2_{[0,T]},
\end{align*}
for $K_i$ defined by induction
\begin{align*}
&K_0= \ max \omega \vee 1,\\
&K_1=1+ \sum_{r=3}^{\infty} \left(\frac{2}{r-2} \right)^{\theta},\\
&K_i=  K_1 \left[ 1 + \sum_{l=1}^i \left( 2K_0^{i/p} K_{k+1-i} + K_i K_{k+1-i } K_0^{\theta} \right) \right].
\end{align*}
\end{theorem}
The following theorem shows that in fact the map $X \rightarrow \hat{X}$ is continuous.
\begin{theorem}
Let $X$ and $Y$ be two almost rough paths of roughness $p$ in
$T^{([p])}(V),$ both of which controlled by a control $\omega,$
that is
\begin{align*}
|X_{s,t}^i|,~~|Y_{s,t}^i| \leq \omega(s,t)^{i/p},~~\forall i=1,...,[p], ~~\forall(s,t) \in \Delta^2_{[0,T]}.
\end{align*}
and for some $\theta >1$ $|(X_{s,t} \otimes X_{t,u})^{i/p}- X^{i/p}_{s,u} | \leq \omega (s,u)^{\theta}$ for all $(s,t),$ $(t,u)\in \Delta^2_{[0,T]},$ $i=1,... ,[p],$ with the same inequality also holding for $Y.$ Suppose that
\begin{align*}
|X^i_{s,t} -Y^i_{s,t} |\leq \varepsilon \omega(s,t)^{i/p},~~~\forall i=1,...,[p],~~~ \forall (s,t) \in \Delta^2_{[0,T]},
\end{align*}
then
\begin{align*}
|\hat{X}^i_{s,t} - \hat{Y}_{s,t}^i | \leq B_i( \varepsilon) \omega(s,t)^{i/p},~~~\forall i=1,...,[p],~~\forall (s,t) \in \Delta^2_{[0,T]};
\end{align*}
where the $B_i$ are defined inductively by

\begin{align*}
&B_1( \varepsilon)= \varepsilon+3 \left\{ \sum_{r=3}^{\infty}
\left[ \varepsilon\left( \frac{2}{r-2}\right)^{1/p} \right] \wedge \left[ \left( \frac{2}{r-2}\right)^{\theta} K_0^{\theta} \right] + \varepsilon \wedge K_0^{\theta} \right\},\\
&B_{k+1}( \varepsilon) = \varepsilon +3 \sum_{r=2}^{\infty} A_k(r, \varepsilon)\end{align*}
and
\begin{align*}
&A_k(r,\varepsilon)= \min
 \{
 K_0^{\theta} (\frac{2}{r-2})^{\theta} ( 1+2 \sum_{i=1}^k K_{k+1-i}K_0^{i/p} + K_i K_{k+1 -i} K_0^{\theta} ),\\
&~~~~  ( \frac{2}{r-2})^{(k+1)/p} [ \sum_{i=1}^k [B_i( \varepsilon) (1+ K_{k+1-i}) + K\varepsilon ] + \varepsilon ] \},\\
&
K=\max_{i=1,...[p]}K_i.
\end{align*}
\end{theorem}
\subsubsection{Spaces of rough path}
Let $C(\Delta^2_{[0,T]}, T^{(n)}(V)) $ denote the set of all continuous functions from the simplex $\Delta^2_{[0,T]}$ into the truncated tensor algebra $(T^{(n)}(V),|.|).$ If $X \in C(\Delta^2_{[0,T]}, T^{(n)}(V)) ,$  then we may write
$$ X_{s,t}=( X_{s,t}^0,...,X_{s,t}^n),~~~\forall (s,t) \in \Delta^2_{[0,T]},$$
where $X^i_{s,t} \in V^{\otimes i}$ is the $i$th component of $X$ (also called the $i$th level path of $X).$
The subset of the functions of $C(\Delta^2_{[0,T]}, T^{(n)}(V)) $ such that $X^0_{s,t}=1$ is denoted $C_0(\Delta^2_{[0,T]}, T^{(n)}(V)).$
\begin{definition}
A function  $X\in C_0(\Delta^2_{[0,T]}, T^{(n)}(V))$ is said to have finite total $p$-variation if
\begin{align*}
\sup_D \sum_l |X^i_{t_{l-1},t_l}|^{p/i} < \infty,~~i=1,...,n,
\end{align*}
where $ sup_D$ runs over all finite subdivisions of $[0,T].$
\end{definition}

\begin{proposition}
Let $p \geq 1$ be a constant, and let $X \in C_0(\Delta^2_{[0,T]}, T^{(n)}(V))$ satisfy Chen's identity (\ref{chenrule}) ( i.e. $X$ is a multiplicative functional in $T^{(n)}(V)$ of order $n).$ If $X$ has a finite $p$-variation, then
\begin{align*}
\omega(s,t)= \sum_{i=1}^n \sup_{D_{[s,t]}} \sum_l |X_{t_{l-1},t_l}|^{p/i},~~\forall(s,t) \in \Delta^2_{[0,T]}
\end{align*}
is a control function, and
\begin{align*}
|X^i_{s,t}|\leq \omega(s,t)^{i/p},~~\forall i=1,...,n,~~~\forall(s,t) \in \Delta^2_{[0,T]}.
\end{align*}
\end{proposition}
Let $C_{0,p}(\Delta^2_{[0,T]}, T^{(n)}(V)) $ denote the subspace of all
$X\in C_0(\Delta^2_{[0,T]}, T^{(n)}(V)) $ with finite $p$ variation. It is
clear that  $C_{0,p}(\Delta^2_{[0,T]}, T^{(n)}(V)) $ is a metric space.
The $p$ variation metric $d_p$ on $C_{0,p}(\Delta^2_{[0,T]},
T^{([p])}(V)))$ is defined by
\begin{align}\label{defdist}
d_p(X,Y)= \max_{i=1,...,[p]}  \sup_D \left( \sum_l |X^i_{t_{l-1},t_l} - Y^i_{t_{l-1},t_l} |^{p/i}\right)^{i/p}.
\end{align}
The space $( \Omega_p(V), d_p)$ is a complete metric space.
\begin{definition}
Geometric rough paths with roughness $p$ are the rough paths in
the closure of smooth rough path under the $p-$ variation
distance.
\end{definition}
The space of all geometric rough paths with roughness $p$ is denoted by $G\Omega_p(V).$
\subsection{Integration theory degree 2}
Let $p \in [2,3[.$
Let $W$ be ${\mathbb R}^n,$ $n \geq 1.$ Let $\alpha ~: V \rightarrow {\mathbf L}(V,W)$ be a function which sends elements of $V$ linearly to $W$-valued one-forms on $V.$ Suppose that $\alpha $ posseses all $k$th continuous derivatives $d^k\alpha$ up to the degree 3 and denote $\alpha^i=d^i \alpha,~~i=1,2,3.$\\
Let $X \in \Omega_{p}(V)$ and let $X_{s,t}=(1,X^1_{s,t},X^2_{s,t}).$ The almost rough path which defines the path integral $\int \alpha(X) dX$ is $Y \in C_{0}(\Delta^2_{[0,T]} ,T^{(2)}(W))$ where $Y_{s,t}=(1,Y^1_{s,t},Y_{s,t}^2)$ and
\begin{align*}
&Y_{s,t}^1 = \alpha^1 (X^1_{0,s}).X^1_{s,t} + \alpha^2(X_{0,s}^1).X^2_{s,t},\\
&Y_{s,t}^2 = \alpha^1 (X^1_{0,s})\otimes \alpha^1 (X^1_{0,s}) .X^2_{s,t}.
\end{align*}
\begin{theorem}
Theorem 5.2.1 and remark 5.3.1 of \cite{LyQi}\\
Let $\alpha ~: V \rightarrow {\mathbf L}(V,W).$  Suppose that
$\alpha $ possesses all   $k$th continuous derivatives $d^k\alpha$
up to the degree 3 and
\begin{align*}
|d^i \alpha(\xi)|_{{\mathbf L}( V \times ...\times V,W)}\leq M( 1+ |\xi|),~~i=1,...,3,~~\forall \xi \in V.
\end{align*}
 Assume that $X\in \Omega_p(V)$ is controlled by $\omega,$
namely
\begin{align*}
|X_{s,t}^i|\leq \omega(s,t)^{i/p},~~i=1,2,~~\forall (s,t) \in \Delta^2_{[0,T]}.
\end{align*}

Then $Y$ is an almost rough path with roughness $p$ in $T^{(2)}(W)$ with control $\omega$ and $\theta=3/p$ i.e.
there exists a universal constant $C$ such that
\begin{align*}
|(Y_{s,t} \otimes Y_{t,u})^i -Y^i_{s,u}|\leq C M \omega(s,u)^{3/p}.
\end{align*}
\end{theorem}
\begin{definition}
Let $X \in \Omega_p(V).$ Then the integral of the one-form $\alpha$ against the rough path $X,$ denoted by $\int \alpha (X) dX,$ is the unique rough path with roughness $p$ in $T^{(2)}(W)$ associated to the almost rough path
$Y \in C_{0}(\Delta^2_{[0,T]} ,T^{(2)}(W))$ where $Y_{s,t}=(1,Y^1_{s,t},Y_{s,t}^2)$ and
\begin{align*}
&Y_{s,t}^1 = \alpha^1 (X^1_{0,s}) .X^1_{s,t} + \alpha^2(X_{0,s}^1).X^2_{s,t},\\
&Y_{s,t}^2 = \alpha^1 (X^1_{0,s})\otimes \alpha^1 (X^1_{0,s}). X^2_{s,t}.
\end{align*}
\end{definition}
\begin{theorem}\label{integration}
Let $\alpha ~: V \rightarrow {\mathbf L}(V,W).$  Suppose that
$\alpha $ possesses all   $k$th continuous derivatives $d^k\alpha$
up to the degree 3 and
\begin{align*}
|d^i \alpha(\xi)|_{{\mathbf L}( V \times ...\times V,W)}\leq M( 1+ |\xi|),~~i=1,...,3,~~\forall \xi \in V.
\end{align*}
 Assume that $X,~\hat{X}\in \Omega_p(V)$ is controlled by $\omega,$
namely
\begin{align*}
|X_{s,t}^i|, ~|\hat{X}_{s,t}^i|\leq \omega(s,t)^{i/p},~~i=1,2,~~\forall (s,t) \in \Delta^2_{[0,T]},
\end{align*}
and
\begin{align*}
|X_{s,t}^i- \hat{X}_{s,t}^i|\leq \varepsilon \omega(s,t)^{i/p},~~i=1,2,~~\forall (s,t) \in \Delta^2_{[0,T]}.
\end{align*}

Then
\begin{align*}
|\int_s^t  \alpha(X_{0,u})d X_u^i -\int_s^t  \alpha(\hat{X}_{0,u})d\hat{ X}_u^i |\leq K \varepsilon  M \omega(s,u)^{i/p},
\end{align*}
for all $(s,t) \in \Delta^2_{[0,T]}$ and $i=1,2,$ where $K$ is a constant
which is  polynomial in $M,$ $\max \omega.$
\end{theorem}
\subsection{It\^o maps : rough path with $2 \leq p < 3$}
\subsubsection{Framework}
Let $V={\mathbb R}^d$ and $W={\mathbb R}^n.$ Let $f~:~V
\rightarrow {\mathbf L}(V,W)$ be a function, which can be viewed
as a map sending vector of $V$ linearly to a vector field on $W.$
Consider the following differential equation ( initial value
problem)
\begin{align}\label{ed-6.1}
 dY_t&=f(Y_t) dX_t,\nonumber\\
Y_0&=y_0.
\end{align}
Since the integral $\int f(Y)dX$ for rough path $X,~Y$ make no sense generally we are not able to iterate differential equation eqn (\ref{ed-6.1}) to obtain the unique solution directely. To overcome this difficulty, the idea is to combine $X$ and $Y$ together as a new path. We view equation (\ref{ed-6.1}) as
\begin{align}\label{ed-6.2}
dX_t&=dX_t,\nonumber\\
dY_t&=f(Y_t)dX_t,~~Y_0=y_0.
\end{align}
The initial condition of $X$ is irrelevant, and therefore we simply take $X_0=0.$ Define $\hat{f}~:V \oplus W \rightarrow {\mathbf L}( V \oplus W; V \oplus W) $ $(V\oplus W$ is the direct sum of $V$ and $W)$ by
\begin{align*}
\hat{f}(x,y)(v,w)= (v, f(y+ y_0).v),~~\forall (x,y) \in V\oplus W,~~~\forall (v,w) \in V\oplus W.
\end{align*}
Then eqn (\ref{ed-6.2}) can be written in the following more appreciating form
\begin{align}\label{eq-6.4}dZ_t= \hat{f}(Z_t)dZ_t.
\end{align}
Given a rough path $X$ in $V,$ we said that a geometric rough path $Z$ in $V \oplus W$ is a solution to (\ref{ed-6.2}) if
\begin{align*}
&\Pi_V(Z)=X,\\
&Z= \int \hat{f}(Z)dZ,
\end{align*}
where $\Pi_V$ is the projector on $T^{([p])}(V).$
\subsubsection{Existence and uniqueness results}
One can summarize the results proved in  pages 149 to 162 of \cite{LyQi}, when the control is $\omega(t,s)=C|t-s|$ and the vector field, $f,$ and its derivatives are bounded in the following way.
\begin{theorem}\label{exist-unicite-ed-petit}
 Let $f \in C^3( W, {\mathbb L}( V,W))$ be a vector field and let $M$ be a constant such that
\begin{align*}
|d^if(\xi)| \leq M,~~~\forall \xi \in W,~~~i=0,1,2,3~;\\
|d^if(\xi) -d^i f(\eta)|\leq M |\xi- \eta|, ~~~\forall \xi,~\eta  \in W,~~~i=0,1,2,3.
\end{align*}

Let $X$ be a rough path in $T^{(2)}(V)$ with roughness
 $2 \leq p < 3$ controlled by   $\omega,$ where 
 for a constant $C_p,$ $\omega (t,s)= C_p|t-s|,~~(s,t) \in \Delta^2_{[0,1]}.$ 

Then, there exist some constants $\kappa,$  and  $\tilde{\theta}$ 
depending only on $M,$ $p$ such that if
\begin{align}\label{def-tps}
T_1= \kappa\frac{1}{1 +C_p^{\tilde{\theta}}}<1
\end{align}
there exists a unique $Z \in \Omega_p( V \oplus W)$ such that $\Pi_V(Z)=X$ and $Z$ satisfies the following integral equation~:
\begin{align*}
Z^i_{s,t}= \int_s^t \hat{f}(Z) dZ^i,~~~i=1,2,~~~\forall (s,t) \in \Delta^2_{[0,T_1]}.
\end{align*}
Moreover, the following estimation holds
\begin{align}\label{control-ci}
|Z^i_{s,t}| \leq (\frac{1}{2} |t-s|)^{i/p},~~~i=1,2,~~~~\forall (s,t) \in \Delta^2_{[0,T_1]}.
\end{align}
\end{theorem}

Finally, we may extend the solution to the whole interval $[0,1].$
\begin{theorem}\label{exist-unicite-ed}
 Let $f \in C^3( W, {\mathbb L}( V,W))$ be a vector field and let $M$ be a constant such that
\begin{align*}
|d^if(\xi)| \leq M,~~~\forall \xi \in W,~~~i=0,1,2,3~;\\
|d^if(\xi) -d^i f(\eta)|\leq M |\xi- \eta|, ~~~\forall \xi,~\eta  \in W,~~~i=0,1,2,3.
\end{align*}
Let $X$ be a rough path in $T^{(2)}(V)$ with roughness $2 \leq p < 3$ controlled by   $\omega,$ where  for a constant $C_p,$ $\omega (t,s)= C_p|t-s|,~~(s,t) \in \Delta^2_{[0,1]}.$ \\
Then,
there exists a unique $Z \in \Omega_p( V \oplus W)$ such that $\Pi_V(Z)=X$ and $Z$ satisfies the following integral equation~:
\begin{align*}
Z^i_{s,t}= \int_s^t \hat{f}(Z) dZ^i,~~~i=1,2,~~~\forall (s,t) \in \Delta^2_{[0,1]}.
\end{align*}
Moreover, there exist some constants $\kappa$ and $\tilde{\theta}$ depending only on $M,$ $p$ such that the following estimation holds
\begin{align}\label{control-ci-bis}
|Z^i_{s,t}| \leq \kappa (1+C_p^{\tilde{\theta}})( |t-s|)^{i/p},~~~i=1,2,~~~~\forall (s,t) \in \Delta^2_{[0,1]}.
\end{align}
\end{theorem}
\begin{proof}

 For example, we may solve the integral equation beyond $T_1$ by replacing the initial condition $y_0$ by $Y_{T_1}= \Pi_W(Z)_{0,T_1}.$ 
Then, the solution is defined up to the time
\begin{align*}
S_2= 2T_1
\end{align*}
Moreover, we have
\begin{align*}
|Y_{S_2}| \leq 2 (\frac{T_1}{2})^{1/p} + |y_0| .
\end{align*}
By an iteration procedure, the solution is defined up to the time 1.
Let $N$ be an integer such that  such that $N\kappa\frac{1}{1 +C_p^{\tilde{\theta}}} >1.$
Then $N$ is bounded by a polynomial in $C_p.$ Then using the Chen
rules, identity (\ref{chenrule}) and several times estimation
(\ref{control-ci}) we obtain
\begin{align*}
|Z^i_{s,t}| \leq \kappa (1+C_p^{\tilde{\theta}}){\mathcal P}(N) (|t-s|)^{i/p},~~~i=1,2,~~~~\forall (s,t) \in \Delta^2_{[0,1]},
\end{align*}
where ${\mathcal P}(N)$ is polynomial in $N$ and then in $C_p.$
\end{proof}

\end{document}